\documentclass[11pt,times]{article}

\usepackage{graphicx,color}

\usepackage[colorlinks=true, breaklinks=true, pdfstartview=FitV, linkcolor=red, citecolor=blue, urlcolor=blue]
{hyperref}

\newcommand{\Vec}[1]{\mbox{\boldmath$#1$}}
\newcommand*\ii{ \mathrm{i\!i}\, }

\begin{document}

\title{Higher Order Approximation to the Hill Problem Dynamics about the Libration Points.\thanks{%
Commun Nonlinear Sci Numer Simulat (2018)
\href{https://doi.org/10.1016/j.cnsns.2017.12.007}{https://doi.org/10.1016/j.cnsns.2017.12.007} (Pre-print version)}}

\author{
Martin Lara\thanks{ETSI Aeron\'auticos, UPM, Madrid, Spain},
Iv\'an L.~P\'erez\thanks{Dep.~of Mechanical Engineering, University of La Rioja, Spain},
Rosario L\'opez\thanks{Center for Biomedical Research of La Rioja, Spain}
}

\date{}

\maketitle

\begin{abstract}
An analytical solution to the Hill problem Hamiltonian expanded about the libration points has been obtained by means of perturbation techniques. In order to compute the higher orders of the perturbation solution that are needed to capture all the relevant periodic orbits originated from the libration points within a reasonable accuracy, the normalization is approached in complex variables. The validity of the solution extends to energy values considerably far away from that of the libration points and, therefore, can be used in the computation of Halo orbits as an alternative to the classical Lindstedt-Poincar\'e approach. Furthermore, the theory correctly predicts the existence of the two-lane bridge of periodic orbits linking the families of planar and vertical Lyapunov orbits. 
\end{abstract}

\section{Introduction}

The restricted three-body problem, in which a body of negligible mass evolves under the gravitational action of two massive bodies whose relative motion is Keplerian (see \cite{Szebehely1967}, for instance), is a useful approximation to the real dynamics of planetary satellites and minor bodies of the solar system, and is customarily used in the first steps of mission designing for some artificial satellites \cite{Gomezetal1993,Gomezetal1998,Koonetal2001}. Of specific interest in astrodynamics is the use of trajectories related to the libration point dynamics \cite{AlessiGomezMasdemont2009,Gomezetal2003,LeiXu2014,Topputo2016}. These points are relative equilibria in a frame rotating with the Keplerian rotation of the primaries.
\par

Rather than pursuing the general solution of the dynamics about these points, the focus has been put on the computation of particular solutions, either numerically, as is the case of the computation of periodic orbits \cite{Henon1965,Henon1965b,BreakwellBrown1979,Howell1984,LaraRussellVillac2007} as well as other invariant manifolds \cite{HowellPernicka1988,GomezMondelo2001}, or analytically, a case in which the Lindstedt-Poincar\'e method is the usual approach \cite{FarquharKamel1973,Richardson1980b,JorbaMasdemont1999,Masdemont2005}. 

In the case of the collinear points, due to the unstable dynamics, a partial reduction to the center manifold is customarily done to remove the hyperbolic directions, which are a consequence of the saddle\,$\times$\,center\,$\times$\,center character of the linear dynamics. After reduction, the center manifold Hamiltonian is of 2-degrees of freedom (DOF) and, therefore, can be approached with the usual tools of nonlinear dynamics \cite{JorbaMasdemont1999,GomezMondelo2001}. 

A common case in the solar system is that the mass ratio between different pairs of celestial bodies (planet-sun, satellite-planet, etc.) is very small. Then, the Hill problem formulation ---a limiting case of the restricted three-body problem in which the two massive bodies evolve in circular orbits about the mutual center of mass, and the body of negligible mass evolves close to the less massive body when compared to the distance between the massive bodies \cite{Szebehely1967}--- releases the dependency of the dynamics on physical parameters. This fact provides a wide generality to this particular formulation, introduced by Hill in his seminal investigations of the Moon's orbit \cite{Hill1878}, whose application encompasses a variety of problems, ranging from the interaction between particles in planetary rings \cite{PetitHenon1986} to the study of the dynamics of planetary satellites \cite{Vashkovyak1999} and about them \cite{LaraRusselVillacl2007Meccanica,LaraPalacianRussell2010}, or the modeling of tethered satellites applications \cite{PelaezLaraetal2012}.

The Hill problem dynamics is studied with the same techniques as the RTBP \cite{Henon1969,Henon1974,Henon2003,Michalodimitrakis1980,ZagourasMarkellos1985,SimoStuchi2000,GomezMarcoteMondelo2005}. Alternatively, it has been recently proposed to study the dynamics about the libration points from the reduction of the Hill problem Hamiltonian to an integrable one \cite{Lara2017}. Indeed, because the libration points are equilibria of the saddle\,$\times$\,center\,$\times$\,center type, the Hamiltonian can be transformed into one such that, after truncation of higher order terms, the hyperbolic part of the quadratic Hamiltonian has been converted to an integral, and the reduced Hamiltonian, which contains the resonant terms of the elliptic part of the quadratic Hamiltonian, is of 1-DOF. The reduced phase space is the sphere \cite{Kummer1976,Cushman1982} and is customarily described in the Hopf variables \cite{Hopf1931}. However, the solutions of the reduced phase space are more insightfully described in Lissajous variables \cite{Deprit1991}, which allow the reconstruction of the periodic solutions by simply evaluating the equilibria solutions for each value of the elliptic anomaly between $0$ and $2\pi$, on the one hand, and ease the computation of the period of the periodic orbit, on the other.
\par

We focus on the Hill problem and extend the normalization of \cite{Lara2017} to higher orders so that it can provide acceptable approximations of the solution for orbits far away from the libration points. 
Celestial mechanicians, as well as physicists, traditionally normalize a system by removing cyclic variables. Thus, the normalization is routinely split into the preliminary reduction to the center manifold and the following removal of short-period effects \cite{CellettiPucaccoStella2015,Lara2017arXiv}. However, the second reduction requires to handle long Fourier series which may become unwieldy at relatively low orders of the perturbation theory. On the other hand, mathematicians have a long experience dealing with perturbed harmonic oscillators in complex variables \cite{Kummer1976,GiorgilliGalgani1978}. With this alternative, the Hamiltonian reduction can be done with a single normalization. Furthermore, normalization in complex variables becomes a simple exercise of polynomial algebra and provides very simple expressions which only require the arithmetic operations that can be carried out directly by the computer's hardware, in this way speeding notably evaluation of the perturbation solution.
\par

The computation of higher orders of the normalization in complex variables allowed us to obtain a single analytical solution which is valid for all the  periodic orbits of the main families of periodic orbits of the Hill problem originated from the libration points. Namely, the families of vertical and planar Lyapunov orbits, the family of Halo orbits, which bifurcates from the family of planar Lyapunov orbits, and, notably, the two lane bridge linking the families of planar and vertical Lyapunov orbits. The later, which exists only for energy values much higher than those of the libration points, is only achieved when the perturbation solution reaches the 6th order, yet corresponding solutions are just rough approximations of the partner periodic orbits. Acceptable approximations to orbits of this family come out only from the 14th order of the perturbation solution.
\par

Computing higher orders in a perturbation approach may be questioned in two ways. On the one hand, the model approached by perturbations is always a simplification of the real dynamics, so the order of the perturbation solution in which the neglected effects of the dynamics would be apparent must be discussed in each particular application. On the other hand, by reasons of simplicity and efficiency, higher orders are customarily approached with the Lie transforms method \cite{Deprit1969,BoccalettiPucacco1998v2} using floating point arithmetic \cite{DepritRom1970}. Because of that, the propagation of the truncation errors in the successive orders of the Taylor series expansion may increase non-negligibly the numerical errors due to the number representation in the computer, a fact that would make nonsense trying to increase the accuracy of the solution by extending the computations beyond a certain order. To mitigate this last issue the perturbation solution is alternatively approached in integer arithmetic. However, while this last approach avoids the accumulation of truncation errors and, therefore, allows to progress exactly, the drawback of using integer arithmetic is the increasing size of the integers to be handled, which grows from order to order of the perturbation solution and may become enormous at relatively moderate orders. As a consequence, the time and memory requirements of the computation of successive orders grow high, thus making the computation of very high orders unpractical. Comparison of the solutions obtained using both techniques helps in estimating the growth of truncation errors of the  floating point arithmetic perturbation solution, and provides a way of extrapolating practical limits for the applicability of a such kinds of solutions.
\par

When the normalization is carried out in floating point arithmetic, the use of complex variables also helps in estimating the accumulation of truncation errors at each order of the perturbation theory. Indeed, when coming back from complex to real variables some residual complex terms will remain in the normalized Hamiltonian, and the size of the higher of the coefficients affecting this residual terms can be taken as an indicator of the truncation errors accumulated in the computations. 
\par

The paper is organized as follows. First, in Sect.~\ref{se:basics}, the Hill problem Hamiltonian is directly derived from the Newtonian dynamics \cite{LaraSanJuan2005}, the origin is translated to a libration point, and the resulting Hamiltonian is expanded in Legendre polynomials in order to present a perturbative arrangement. Next, the linearized dynamics of the problem is discussed in Sect.~\ref{se:linear}; while this part is mostly borrowed from \cite{Lara2017} additional details are given to show that the linear transformation that decouples the linearized dynamics is not unique; still, we adhere to tradition and use the usual transformation \cite{JorbaMasdemont1999} in our computations. It follows the description of the normal form computation in Sect.~\ref{se:normalization}, which is computed exactly using integer arithmetic up to the order 11, and approximately using floating point arithmetic up to the order 20; this section provides estimates of the accumulation of the truncation errors in the floating point arithmetic case due to the physical size of the computer's registers. Also in this section, the equilibria of the reduced phase space are briefly discussed to show their correspondence with well known periodic solutions of the Hill problem. Finally, a variety of tests are presented in Sect.~\ref{se:performance} to illustrate the performance, as well as the limits, of the analytical solution.

\section{Hill problem Hamiltonian about the libration points} \label{se:basics}

Let $P$ and $S$ be two massive bodies, of masses $M$ and $m$, respectively, which, under their mutual gravitational attraction, are evolving in circular orbits about the system's center of mass with constant angular velocity $\Vec{\omega}$. Then, the distance $d$ between $P$ and $S$ remains constant. Let $O$ be a massless body evolving under the gravitational actions of $S$ and $P$, and let $\Vec{R}$ define the position of $O$ with respect to the center of mass of the system. Then, from Newton's gravitational law,
\begin{equation} \label{ainertial}
\frac{\mathrm{d}^2\Vec{R}}{\mathrm{d}t^2}=-\frac{GM}{s^3}\Vec{s}-\frac{Gm}{r^3}\Vec{r}
\end{equation}
where $t$ is the usual time, $\Vec{s}$ is the vector joining $P$ with $O$, of modulus $s$, $\Vec{r}$ is the vector joining $S$, with $O$, $r=\|\Vec{r}\|$, and $G$ is the gravitational constant.

We study the motion of $O$ relative to $S$ in a rotating frame $(S,\Vec{i},\Vec{j},\Vec{k})$ with rotation rate $\Vec{\omega}=\omega\Vec{k}$ defining the $z$ axis direction $\Vec{k}$, the $x$ axis direction $\Vec{i}$ is defined by the direction from $P$ to $S$, and the $y$ axis direction $\Vec{j}$ completes a direct orthogonal frame. Then, $\Vec{R}=q\Vec{i}+\Vec{r}$, where, from the definition of the center of mass
\[
q=\frac{M}{M+m}d,
\]
and, from the derivative of a vector in a rotating frame,
\begin{equation} \label{Darboux}
\frac{\mathrm{d}^2\Vec{R}}{\mathrm{d}t^2}=\ddot{\Vec{r}}+2\omega\Vec{k}\times\dot{\Vec{r}}+\omega^2\Vec{k}\times(\Vec{k}\times\Vec{r})-\omega^2q\Vec{i},
\end{equation}
where dots over vectors mean differentiation in the rotating frame. Therefore, from Eqs.~(\ref{ainertial}) and (\ref{Darboux}),
\begin{equation} \label{arotating}
\ddot{\Vec{r}}+2\omega\Vec{k}\times\dot{\Vec{r}}=\omega^2q\Vec{i}-\omega^2\Vec{k}\times(\Vec{k}\times\Vec{r})
-\frac{GM}{s^3}\Vec{s}-\frac{Gm}{r^3}\Vec{r},
\end{equation}
where
\begin{equation}
\Vec{s}=d\Vec{i}+\Vec{r}=(d+x)\Vec{i}+y\Vec{j}+z\Vec{k},
\end{equation}
and, due to the circular motion of $S$, the centripetal acceleration is 
\begin{equation} \label{centripetal}
\omega^2q=\frac{GM}{d^2}.
\end{equation}

Now, with Hill, we assume that $r\ll{d}$ and $m\ll{M}$. More precisely, we assume that $r/d=\mathcal{O}(\epsilon)$ and $m/M=\mathcal{O}(\epsilon^2)$. Then, $q=d-\mathcal{O}\left(\epsilon^2\right)$ and
\[
\frac{GM}{s^3}\Vec{s}=\omega^2q\frac{\Vec{s}/d}{(s/d)^3}=
\omega^2\left[\left(d-2x\right)\Vec{i}+y\Vec{j}+z\Vec{k}\right]+\mathcal{O}\left(\epsilon^2\right).
\]
Hence, after neglecting higher order terms, Eq.~(\ref{arotating}) is rewritten as the differential equation of the Hill problem
\begin{equation} \label{Hillflow}
\ddot{\Vec{r}}+2\omega\Vec{k}\times\dot{\Vec{r}}=\omega^2(3x\Vec{i}-z\Vec{k})-\frac{Gm}{r^3}\Vec{r}.
\end{equation}
The dot product of $\dot{\Vec{r}}$ and Eq.~(\ref{Hillflow}) can be integrated, to yield
\[
\mathcal{E}=\frac{1}{2}\left(\dot{x}^2+\dot{y}^2+\dot{z}^2\right)-\frac{Gm}{r}-\frac{1}{2}\omega^2\left(3x^2-z^2\right).
\]

After scaling units of length by $(Gm)^{1/3}$ and time by $1/\omega$, Eq.~(\ref{Hillflow}) reads
\begin{eqnarray} \label{ax}
\ddot{x}-2\dot{y} &=& -\frac{x}{r^3}+3x, \\ \label{ay}
\ddot{y}+2\dot{x} &=& -\frac{y}{r^3}, \\ \label{az}
\ddot{z} &=& -\frac{z}{r^3}-z,
\end{eqnarray}
revealing that the Hill problem does not depend on any parameter. Besides, it is simple to check that two equilibria, the so-called libration points, exist at the positions $y=z=0$, $x=\pm\xi$, where $\xi=3^{-1/3}\approx0.69$.

The Hill problem accepts Hamiltonian formulation. Indeed, the flow in Eqs.~(\ref{ax})--(\ref{az}) can be derived from the Hamiltonian
\begin{equation} \label{Jota}
\mathcal{J}=
\frac{1}{2}\left(X^2+Y^2+Z^2\right)+Xy-xY-\frac{1}{r}+\frac{1}{2}\left(r^2-3x^2\right),
\end{equation}
where $X=\dot{x}-y$, $Y=\dot{y}+x$, and $Z=\dot{z}$ are the conjugate momenta to $x$, $y$, and $z$, respectively.
\par

Because of the known symmetries of the Hill problem (see \cite{GomezMarcoteMondelo2005}, for instance), it is enough to study the dynamics about one of the libration points, say $x=\xi$. To do that, the origin is translated to $\xi$. Since a translation is a canonical transformation, after neglecting constant terms the transformed Hamiltonian reads
\begin{equation} \label{H1ham}
\mathcal{H} =\frac{1}{2}\left(X'^2+Y'^2+Z'^2\right)-\left(x'Y'-y'X'\right)+\frac{1}{2}\left(y'^2+z'^2\right)-x'^2-\frac{x'}{\xi^2}-\frac{1}{r},\end{equation}
where
\begin{equation} \label{translation}
x=x'+\xi, \quad y=y', \quad z=z', \quad X=X', \quad Y=Y'+\xi, \quad Z=Z'.
\end{equation}
\par

Now,
\[
r=\sqrt{\left(x'+\xi\right)^2+y'^2+z'^2}=\xi\sqrt{1-2(r'/\xi)\cos\psi+(r'/\xi)^2},
\]
where $\cos\psi=-x'/r'$ and $r'$ is the distance from $O$ to the libration point. Then, for values $r'/\xi<1$, the term $1/r$ can be replaced by the usual expansion in Legendre polynomials $P_{n}(\cos\psi)$, yielding\begin{equation} \label{Hillps}
\mathcal{H} =H_0-\frac{1}{\xi}\sum_{n>0}\left(\frac{r'}{\xi}\right)^{n+2}P_{n+2}\left(\cos\psi\right),
\end{equation}
in which
\begin{equation} \label{H0}
H_0=\mbox{$\frac{1}{2}$}(X'^2+Y'^2)-(x'Y'-X'y')+2(y'^2-2x'^2)+\mbox{$\frac{1}{2}$}(Z'^2+4z'^2),
\end{equation}
whereas the other terms of the Hamiltonian comprise monomials of the form
\begin{equation} \label{monomials}
M_k=Q_k{x'}^{m_1}{X'}^{m_2}{y'}^{m_3}{Y'}^{m_4}{z'}^{m_5}Z'^{m_6}, \quad k=(m_1,m_2,m_3,m_4,m_5,m_6),
\end{equation}
where $Q_k$ are numeric coefficients and $m_i$ ($i=1,\dots,6$) are non-negative integers.
\par

\section{Linear dynamics about the libration points} \label{se:linear}

For small displacements about the libration point we can neglect terms of higher order than $(r'/\xi)^2$, and hence the zeroth order term (\ref{H0}) of the Hamiltonian (\ref{Hillps}) is representative of the dynamics. The last term in Eq.~(\ref{H0}) has the form of a simple harmonic oscillator with frequency $\nu=2$. Therefore, in the linear approximation, the motion in $z$ and $Z$ decouples from the rest of the flow and is made of harmonic oscillations. That is, the equilibria of the Hill problem are of the center type relative to the $z$ direction.
\par

On the other hand, the coupled $x$-$y$ motion results from the integration of a linear differential system with constant coefficients. Indeed, from Hamilton equations,
\begin{equation} \label{flow1}
\left(\dot{x'}, \dot{y'}, \dot{X'}, \dot{Y'}\right)^\tau
=M_1\left(x',y',X',Y'\right)^\tau,
\end{equation}
where $\tau$ means transposition, and
\begin{equation} \label{M1}
M_1= \left(\begin{array}{rrrr} 
 0 & 1 & 1 & 0 \\ 
-1 & 0 & 0 & 1 \\ 
 8 & 0 & 0 & 1 \\ 
 0 &-4 &-1 & 0 
\end{array}\right).
\end{equation}

The general solution of Eq.~(\ref{flow1}) is
\begin{equation} \label{xi}
\left(x',y',X',Y'\right)^\tau=B\exp(\Vec{\lambda}t),
\end{equation}
where $B=(b_{i,j})$ is a 4 $\times$ 4 matrix of arbitrary coefficients, and $\Vec{\lambda}$ is the vector of characteristic exponents $\lambda_j$, $j=1,\dots,4$, which are the eigenvalues of $M_1$. It is simple to check that
\begin{equation} \label{lambda}
\lambda_{1,2}=\pm\lambda, \qquad \lambda=\sqrt{2\sqrt{7}+1},
\end{equation}
are real numbers, thus giving place to an hyperbolic component of Eq.~(\ref{xi}), a saddle direction, whereas
\begin{equation} \label{omega}
\lambda_{3,4}=\pm\omega\ii, \qquad \omega=\sqrt{2\sqrt{7}-1},
\end{equation}
with $\ii=\sqrt{-1}$, are pure imaginary numbers, resulting in an elliptic or center-type component of Eq.~(\ref{xi}). These exponents define the well-known saddle\,$\times$\,center\,$\times$ center type of the libration points of the Hill problem.
\par

Note that the quadratic Hamiltonian
\begin{equation} \label{separable}
K_0=\lambda{x}_1X_1+\frac{1}{2}\left(Y_1^2+\omega^2y_1^2\right)+\frac{1}{2}\left(Z_1^2+\nu^2z_1^2\right),
\end{equation}
is in separate variables and enjoys the same dynamical behavior as Eq.~(\ref{H0}). The Hamiltonian flow stemming from Eq.~(\ref{separable}) is
\begin{equation} \label{flow2}
\left(\dot{x}_1,\dot{y}_1,\dot{X}_1,\dot{Y}_1\right)^\tau = M_2 \left(x_1,y_1,X_1,Y_1\right)^\tau,
\end{equation}
with
\begin{equation} \label{M2}
M_2=\left(\begin{array}{rcrr} 
\lambda &        0 &       0 & 0 \\ 
      0 &        0 &       0 & 1 \\ 
      0 &        0 &-\lambda & 0 \\ 
      0 &-\omega^2\; &       0 & 0
\end{array}\right).
\end{equation}
So it emerges the question if a canonical transformation can be found such that it transforms $H_0$ into $K_0$. The answer is in the affirmative, and, because the equations of motion are linear, the linear transformation
\begin{equation} \label{x1tox}
\left(x',y',X',Y'\right)^\tau =A\left(x_1,y_1,X_1,Y_1\right)^\tau,
\end{equation}
can be computed by solving $A=A(a_{i,j})$ from the underdetermined linear system
\begin{equation} \label{M1AAM2}
M_1A=AM_2,
\end{equation}
which is obtained by equating the right sides of Eqs.~(\ref{flow1}) and (\ref{flow2}), and replacing Eq.~(\ref{x1tox}).
\par

The solution of Eq.~(\ref{M1AAM2}) is expressed as a function of 4 arbitrary coefficients, say
\begin{equation} \label{matrixA}
A=\left(
\begin{array}{cccc}
 -\frac{1}{4} \left(\omega ^2+7\right) a_{4,1} & -\frac{1}{4} \left(\omega
   ^2+3\right) a_{3,4} & -\frac{1}{4} \left(\omega ^2+7\right) a_{4,3} &
   \frac{1}{4} \left(\omega ^2-5\right) a_{4,4} \\ [0.5ex]
 \frac{1}{12} \lambda  \left(\omega ^2+3\right) a_{4,1} & -\frac{1}{4}
   \left(\omega ^2-9\right) a_{4,4} & -\frac{1}{12} \lambda  \left(\omega
   ^2+3\right) a_{4,3} & -\frac{1}{4} \left(\omega ^2+7\right) a_{3,4} \\ [0.5ex]
 -\frac{1}{3} \lambda  \left(\omega ^2+6\right) a_{4,1} & \left(2 \omega
   ^2-9\right) a_{4,4} & \frac{1}{3} \lambda  \left(\omega ^2+6\right) a_{4,3}
   & a_{3,4} \\ [0.5ex]
 a_{4,1} & \left(\omega ^2+6\right) a_{3,4} & a_{4,3} & a_{4,4} \\[0.5ex]
\end{array}
\right)
\end{equation}
where
\begin{equation} \label{detA}
|A|= -\frac{14}{3} \lambda  a_{4,1} a_{4,3} \left[(10 \omega ^2+63)a_{3,4}^2+\omega ^2 a_{4,4}^2\right]\ne0,
\end{equation}
for an invertible transformation, and hence $a_{4,1}\ne0$, $a_{4,3}\ne0$, while $a_{3,4}$ and $a_{4,4}$ cannot vanish at the same time .
\par

Additionally, the transformation in Eq.~(\ref{x1tox}) must be canonical. Because of the linear character of the transformation, it happens that its Jacobian is also $A$. Therefore, the canonicity is expressed as $AJA^\tau=J$, where $J$ is the symplectic matrix of dimension 4, yielding the two additional constraints 
\begin{equation} \label{constraint}
\lambda a_{4,1}a_{4,3}=-\frac{23-5 \omega ^2}{28}, \qquad
a_{3,4}^2+(7\omega^2-30)a_{4,4}^2=\frac{23-5\omega^2}{42},
\end{equation}
which make $|A|=1$ in Eq.~(\ref{detA}), and left undetermined only two coefficients, say $a_{4,1}$ and $a_{4,4}$. 
\par

Among the different possibilities, and in view of the similitudes of columns 1 and 3 of the matrix $A$, it seems natural to choose $a_{3,4}=0$, $a_{4,3}=-a_{4,1}$. In view of Eq.~(\ref{constraint}), it yields
\[
a_{4,1}=\frac{\lambda}{\sigma}(\lambda^2-7), \qquad
a_{4,4}=-\frac{1}{\tau}(\omega^2+7),
\]
with
\[
\sigma =4 \sqrt{\lambda  \left(2 \lambda ^2-9\right)}, \qquad \tau =2\sqrt{2(2 \omega ^2+9)},
\]
which leads to the usual transformation matrix
\begin{equation} \label{Agerard}
A=\left(
\begin{array}{cccc}
 2 \lambda/\sigma & 0 & -2 \lambda/\sigma & 2/\tau \\ 
 (\lambda^2-9)/\sigma & -(\omega^2+9)/\tau & (\lambda^2-9)/\sigma & 0 \\
 (\lambda^2+9)/\sigma & (9-\omega^2)/\tau & (\lambda^2+9)/\sigma & 0 \\
 \lambda(\lambda^2-7)/\sigma & 0 & \lambda(7-\lambda^2)/\sigma & -(\omega^2+7)/\tau
\end{array}
\right),
\end{equation}
based on the eigenvector decomposition of $M_1$ in Eq.~(\ref{M1}), cf.~\cite{JorbaMasdemont1999} (see, also, \cite{DepritDelie1965}).
\par

The discussion of the particular merits of the different transformations derived from the bi-parametric family defined by Eqs.~(\ref{matrixA}) and (\ref{constraint}) is not made here. However to further stress that the use of Eq.~(\ref{Agerard}) is not a requirement to achieve the reduction of the quadratic Hamiltonian (\ref{H0}) to the separable form in Eq.~(\ref{separable}), we write explicitly an alternative choice. Namely,
\[
a_{4,1}=-\frac{23-5 \omega ^2}{112\lambda(\omega^2-4)}, \qquad
a_{4,4}=0,
\]
that yields
\begin{equation} \label{decoupling}
A=
\left(
\begin{array}{c@{\quad}c@{\quad}c@{\quad}c}
 \frac{1}{6048}\lambda(11\omega^2+81) &-\frac{1}{2}\gamma(\omega^2+3) &-(\omega^2-1) & 0 \\[1ex]
-\frac{1}{672}(\omega^2+15) & 0 & -\lambda  \left(5-\omega^2\right) &-\frac{1}{2}\gamma(\omega^2+7) \\[1ex]
 \frac{1}{336} \left(5 \omega ^2+33\right) & 0 & 4 \lambda  & 2 \gamma \\[1ex]
-\frac{1}{3024}\lambda(27-\omega^2) & 2 \gamma  \left(\omega^2+6\right) & 4 \left(\omega ^2-4\right) & 0 \\
\end{array}
\right),
\end{equation}
where $\gamma=\sqrt{(23-5\omega)/168}$.

\section{Normalization} \label{se:normalization}

The instability of the libration points due to their saddle component can be skipped by choosing suitable initial conditions. Indeed, due to the center\,$\times$\,center part, when the hyperbolic part $J=x_1X_1$ of the quadratic Hamiltonian (\ref{separable}) vanishes the resulting motion will be periodic or quasi-periodic. This dynamics is not constrained to the linear approximation and can be extended to the whole Hamiltonian
\begin{equation} \label{Kam}
\mathcal{K}=\sum_{n\ge0}\frac{1}{n!}K_n(x_1,y_1,z_1,X_1,Y_1,Z_1),
\end{equation}
which is obtained after applying the transformation (\ref{x1tox}) to Hamiltonian (\ref{Hillps}). The procedure for removing the hyperbolic components of Eq.~(\ref{Kam}) is called the reduction to the center manifold, and consist of a partial normalization that, after truncation to some order, converts $J$ into an integral, which is followed by constraining the motion to the manifold $J=0$ \cite{GomezJorbaMasdemontSimo1991}.
\par

A following removal of the short-period terms normalizes the Hamiltonian \cite{Lara2017}. However, the normalization of Eq.~(\ref{Kam}) does not need to be split into two different parts and, on the contrary, can be achieved with a single transformation. This is the approach we take here.
\par

First of all, in view of the frequencies $\nu$ and $\omega$ are very close, we introduce a detuning parameter \cite{Henrard1970}
\begin{equation} \label{delta}
\delta=1-(\nu/\omega)^2=\frac{23-8\sqrt{7}}{27}\approx\frac{1}{10}\xi.
\end{equation}
Then, the zeroth order term (\ref{separable}) is split into $K_0=K_0^*-\frac{1}{2}\delta\omega^2z_1^2$, and we rearrange the Hamiltonian (\ref{Kam}) as 
\begin{equation} \label{KK}
\mathcal{K}=\sum_{n\ge0}\frac{1}{n!}K_n^*,
\end{equation}
with
\begin{eqnarray} \label{detuned}
K_0^* &=& \lambda{x}_1X_1+\frac{1}{2}\left[Y_1^2+Z_1^2+\omega^2\left(y_1^2+z_1^2\right)\right], \\ \label{detuned1}
K_1^* &=& K_1-\frac{1}{2}\delta\omega^2z_1^2, \\
K_n^* &=& K_n \quad (n>1),
\end{eqnarray}
where the $K_n^*$ comprise monomials in the subindex $1$ variables of the type given in Eq.~(\ref{monomials}). Note that, because of the detuning, $K_1^*$ is no longer an homogeneous polynomial of degree 3. With this artifact, the center components of Eq.~(\ref{detuned}) have the form of an elliptic oscillator of frequency $\omega$.
\par

The normalization of the Hamiltonian (\ref{KK}) is more easily achieved in complex variables \cite{DepritDelie1965,Kummer1976,GiorgilliGalgani1978}. Thus, the transformation
\begin{equation} \label{tocomplex} \begin{array}{rclrcl}
x_1 &=& u', &\qquad X_1 &=& U', \\[0.5ex]
y_1 &=& \displaystyle\frac{1}{\sqrt{2\omega}}\left(v'+\ii V'\right), &\qquad Y_1 &=& \displaystyle\sqrt{\frac{\omega}{2}}\left(V'+\ii v'\right), \\[2ex]
z_1 &=& \displaystyle\frac{1}{\sqrt{2\omega}}\left(w'+\ii W'\right), &\qquad Z_1 &=& \displaystyle\sqrt{\frac{\omega}{2}}\left(W'+\ii w'\right),
\end{array}\end{equation}
with $\ii=\sqrt{-1}$, is applied first. Hence,
\begin{equation} \label{complex}
K_0^*=\lambda{u}'U'+\ii\omega\left(v'V'+w'W'\right),
\end{equation}
while the remainder terms of the Hamiltonian stay as monomials of the type of Eq.~(\ref{monomials}) in the new real-complex variables. The form of the zeroth order Hamiltonian (\ref{complex}) converts the normalization process in an elementary problem of polynomial algebra (see \cite{Lara2017} for details).
\par

\subsection{The reduced dynamics} \label{ss:reduced}

The normalization is implemented using the Lie transforms procedure \cite{Deprit1969,BoccalettiPucacco1998v2}, and yields a Hamiltonian in new variables $(u,v,w,U,V,W)$ in which, after truncation to some desired order, the hyperbolic and elliptic components of the quadratic Hamiltonian, viz. 
\begin{equation} \label{integrals}
J=uU, \qquad L=\ii(vV+wW),
\end{equation}
become formal integrals. Therefore, the normalized Hamiltonian is of one degree of freedom, and hence integrable. Furthermore, to skip the hyperbolic instability, we choose initial conditions $u=U=0$ in this way constraining the dynamics to the center manifold $J=0$. In consequence, the final, normalized Hamiltonian is of the form
\begin{equation} \label{normalized}
\mathcal{N}=\omega{L}+\sum_{n\ge1}\frac{1}{n!}N_n(v,w,V,W),
\end{equation}
where terms $N_n$ are polynomials in the complex variables, which are no longer homogeneous because of the detuning made in Eqs.~(\ref{detuned})--(\ref{detuned1}). The first terms of the normalized Hamiltonian (\ref{normalized}) are printed below:
\begin{eqnarray} \label{N1}
N_1 &=& -\frac{1}{2}\ii\omega\delta w W, 
\\ \label{N2}
N_2 &=& 3 \xi\frac{1561 \omega ^2+701}{38696}(vV)^2 + \xi\frac{15235\omega^2-41269}{38696}\left[(vW)^2+(Vw)^2\right] \\ \nonumber
&& +\xi\frac{89-23\omega^2}{42}vVwW + \xi\frac{484\omega^2-193}{4146}(wW)^2 - \frac{1}{4}\ii\omega\delta(\delta+2)wW, 
\\ \label{N3}
N_3 &=& \delta\xi\frac{12286157 \omega ^2-60125147}{13369468}\left[(vW)^2+(Vw)^2\right] -\frac{3}{8}\ii\omega\delta^3wW \\ \nonumber
&& +\delta\xi\frac{619-109 \omega^2}{84}vVwW + 2\delta\xi\frac{257575 \omega^2-737704}{1432443}(wW)^2,
\\ \label{N4}
N_4 &=& 13\delta^2\xi\frac{52000002667\omega^2-333151697455}{166289442984}\left[(vW)^2+(Vw)^2\right] \\ \nonumber
&& +\delta^2\xi\frac{10147-1501\omega^2}{252}vVwW + 1184\delta^2\xi\frac{4042054 \omega^2-20055643}{2969454339}(wW)^2 \\ \nonumber
&& -\frac{15}{16}\ii\omega\delta^4wW 
+\ii\xi^2\omega\frac{21541881606067 \omega^2+2841023083259}{51605157139368}(vV)^3 \\ \nonumber
&& +\ii\xi^2\omega\frac{50988038481433\omega^2-196099405389331}{10476234907992}w W\left[(vW)^2+(Vw)^2\right] \\ \nonumber
&& -\ii\xi^2\omega\frac{214918405892794\omega^2-1084821213782455}{66349487750616}v V\left[(vW)^2+(Vw)^2\right] \\ \nonumber
&& +5\ii\xi^2\omega\frac{49355491209497 \omega^2-212878468773272}{31428704723976}(vV)^2 w W \\ \nonumber
&& -\ii\xi^2\omega\frac{4997434564153693 \omega^2-21286323197607436}{597145389755544}vV(wW)^2 \\ \nonumber
&& +\ii\xi^2\omega\frac{71823531673 \omega ^2-119333443606}{160350534306}(wW)^3.
\end{eqnarray}
\par

Note that the integral $L$ is not easily identified in the summands $N_n$ of the normalized Hamiltonian (\ref{normalized}), which seems to remain as a 2 degrees of freedom Hamiltonian in the complex variables. The use of Hopf variables \cite{Hopf1931}, given by the transformation
\begin{equation} \label{hopfComplex}
I_1 = \frac{\ii}{2}(wW-vV), \quad
I_2 =-\frac{\ii}{2}(vW+wV), \quad
I_3 = \frac{1}{2}(vW-wV),
\end{equation}
with the constraint
\begin{equation} \label{sphere}
I_0^2=I_1^2+I_2^2+I_3^2=\mbox{$\frac{1}{4}$}L^2,
\end{equation}
definitely helps in disclosing the formal integral $L=2I_0$, as well as in describing the reduced phase space, which is the sphere \cite{Kummer1976,Cushman1982}. Indeed, Eqs.~(\ref{N1})--(\ref{N4}) are trivially expressed in Hopf variables by using the relations
\begin{equation} \label{sphere}
wW=-\ii(I_1+I_0), \quad vV=\ii(I_1-I_0), \quad vW=\ii{I}_2+I_3, \quad Vw=\ii{I}_2-I_3.
\end{equation}
\par

Then, Eq.~(\ref{normalized}) takes the form
\[
\mathcal{N}=\omega{L}+\sum_{n\ge1}\frac{1}{n!}N_n(I_1,I_2,I_3;L),
\]
from which we derive the Hamiltonian flow $\dot{I}_i=\{I_i;\mathcal{H}\}$, $i=1,2,3$, where curly brackets represent the Poisson bracket operator. We find
\begin{eqnarray} \label{I1p}
\dot{I}_1 &=& I_2I_3\sum_nF_{1,n}(I_1,I_2,I_3;L), \\ \label{I2p}
\dot{I}_2 &=& I_3\sum_nF_{2,n}(I_1,I_2,I_3;L), \\ \label{I3p}
\dot{I}_3 &=& I_2\sum_nF_{3,n}(I_1,I_2,I_3;L).
\end{eqnarray}
Hence, points
\begin{equation} \label{vyp}
\left(\pm\mbox{$\frac{1}{2}$}L,0,0\right),
\end{equation}
on the sphere are always equilibria; the plus sign corresponds to vertical Lyapunov orbits and the minus sign to planar Lyapunov orbits. Besides, those points
\begin{equation} \label{halo}
(I_1,0,\pm{I}_3),
\end{equation}
such that Eq.~(\ref{I2p}) vanish, that is
\begin{equation} \label{haloRoot}
\sum_nF_{2,n}(I_1,0,I_3;L)=0,
\end{equation}
are also equilibria. This new equilibria stem from the point $(-\frac{1}{2}L,0,0)$ in a pitchfork bifurcation at the value $L=L_\mathrm{h}$ given by the root $\sum_nF_{2,n}(\frac{1}{2}L,0,0;L)=0$, and correspond to Halo orbits. Finally, points
\begin{equation} \label{bridge}
(I_1,\pm{I}_2,0),
\end{equation}
on the sphere such that Eq.~(\ref{I3p}) vanish, that is
\begin{equation} \label{bridgeRoot}
\sum_nF_{3,n}(I_1,I_2,0;L)=0,
\end{equation}
are equilibria as well. Computation of the roots $L_\mathrm{b1}$ and $L_\mathrm{b2}$ of the equation $\sum_nF_{3,n}(\frac{1}{2}L,0,0;L)=0$ shows that they stem from $(-\frac{1}{2}L_\mathrm{b2},0,0)$ in a pitchfork bifurcation, and collapse into $(\frac{1}{2}L_\mathrm{b1},0,0)$. These equilibria on the sphere correspond to the two-lane bridge of periodic orbits linking planar and vertical Lyapunov orbits. Interested readers are referred to \cite{Lara2017} for full details on the discussion of the reduced phase space as well as basic references on the topic.

\subsection{Computational issues}

The normalization is computed exactly by avoiding decimal expansions of the involved rational and irrational numbers. The irrational numbers $\xi$, $\lambda$, $\omega$ and $\delta$ are handled formally, an their respective powers are simplified as mach as possible. However, the size of the integer numbers involved in the rational coefficients of the monomials grows from order to order, as can be observed in Eqs.~(\ref{N1})--(\ref{N4}), soon causing memory allocation to become a serious issue, with the consequent rapid increase of computing time. In this way, we only succeeded in extending the computations to the eleventh order, in which the rational coefficients may involve integer numbers of more than 100 digits. We hasten to say that we relied on commercial, general purpose, symbolic algebra tools in our computations; development of specific manipulators by experts could, of course, ease considerably the task \cite{Henrard1988,SanJuan1994,Jorba1999}.
\par

\begin{figure}[htbp]
\centering \includegraphics[scale=0.9, angle=0]{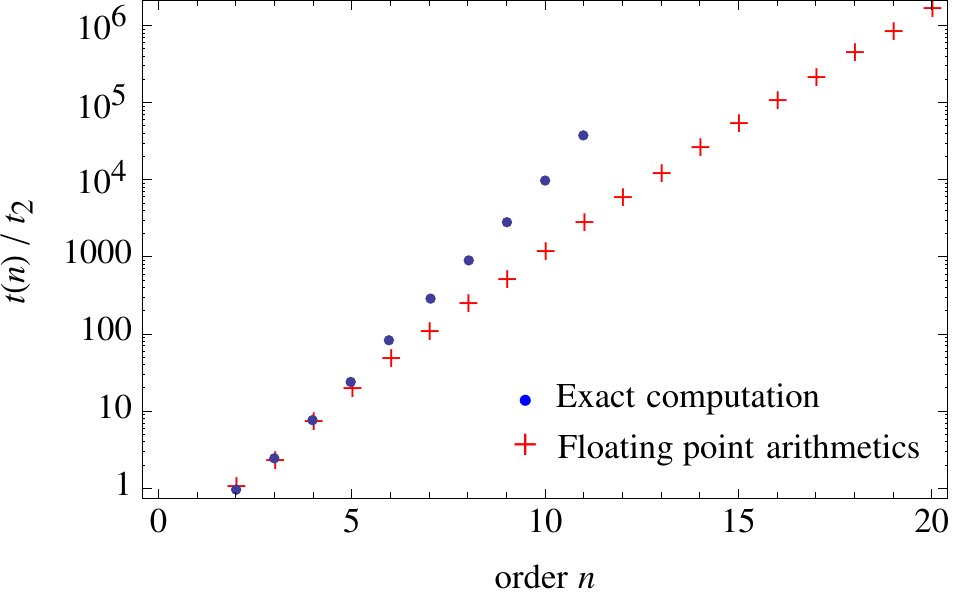}
\caption{Computation time $t(n)$ of the different orders of the perturbation solution scaled by the time $t_2$ spent in the computation of the second order solution.}
\label{f:trationalreal}
\end{figure}

On the other hand, the use of floating point arithmetic expedites computations notably, but at the cost of introducing truncation errors due to the physical length of the computer's registers. The computation time still grows exponentially with the order of the theory, but at a lower rate. This fact is illustrated in Fig.~\ref{f:trationalreal}, where it is shown that when using floating point arithmetic the computation time $t$ grows with the order $n$ roughly as $t(n)\approx0.336\exp(0.788n)t_2$, $n>2$, where $t_2$ is the time spent into the computation of the second order terms of the normalized Hamiltonian and generating function, whereas in the case of exact computations using integer arithmetic it grows as $t(n)\approx0.073\exp(1.18n)t_2$, $n>2$. In fact, the time employed in the exact computation of the order 11 of the perturbation solution was almost $40\,000$ times longer than $t_2$, whereas $t(11)$ was only $2\,500t_2$ in the floating point case; that is, approximately 16 times faster.
\par

\begin{figure}[htbp]
\centering \includegraphics[scale=0.9, angle=0]{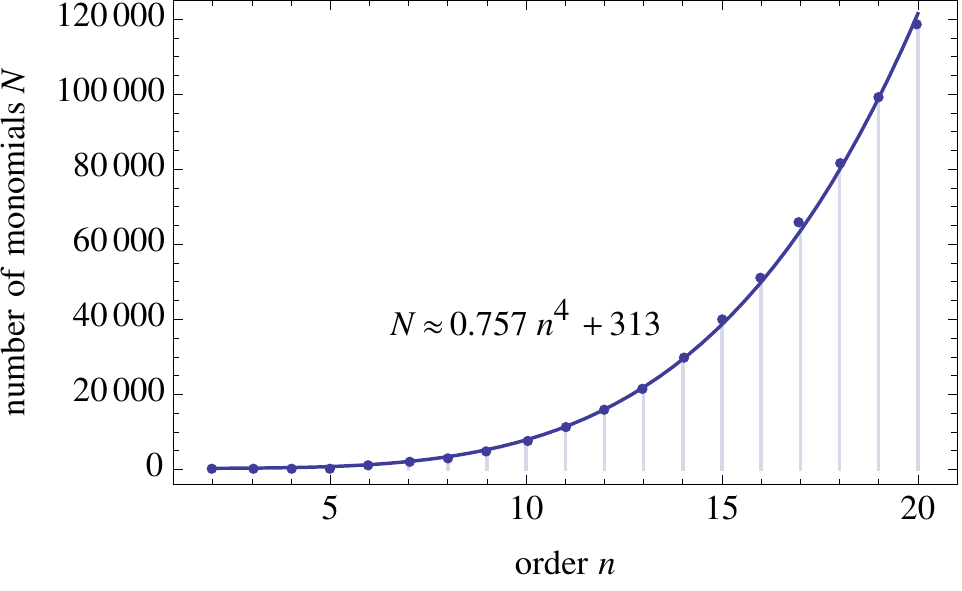}
\caption{Number of monomials involved in each different order of the generating function.}
\label{f:monomialsW}
\end{figure}

Memory handling issues are definitely less severe with the floating point approach. However, the number of terms to be evaluated by consecutive higher orders of the perturbation solution grows roughly with the quartic power of the order, and soon becomes enormous, as shown in Fig.~\ref{f:monomialsW}. Hence, we did not progress in our computations further than the order 20. 
\par

\begin{figure}[htb]
\centering \includegraphics[scale=0.9, angle=0]{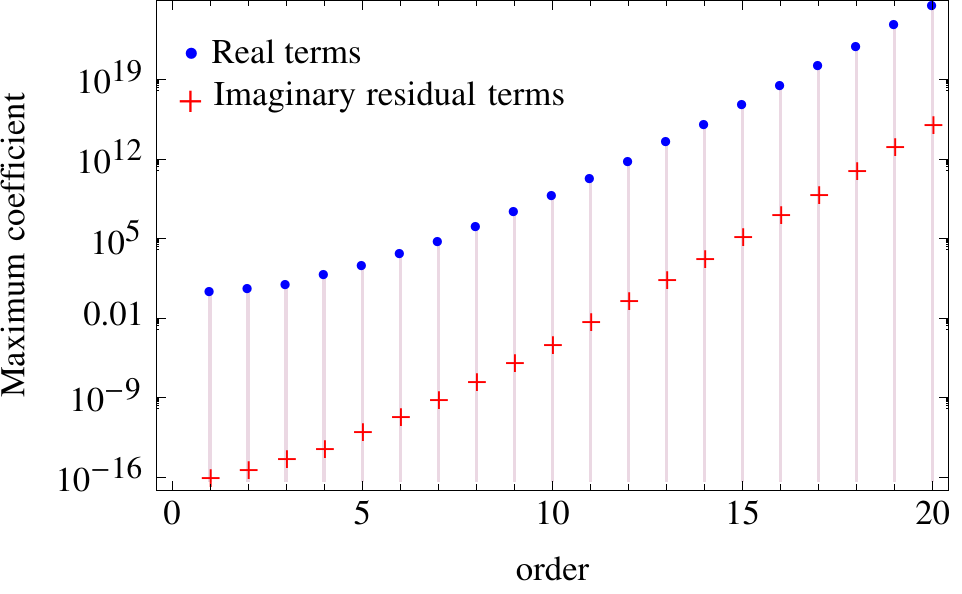}
\caption{Evolution of the maximum coefficients of the generating function with the order of the perturbation theory.}
\label{f:maxCoefW}
\end{figure}

\begin{figure}[htb]
\centering \includegraphics[scale=0.9, angle=0]{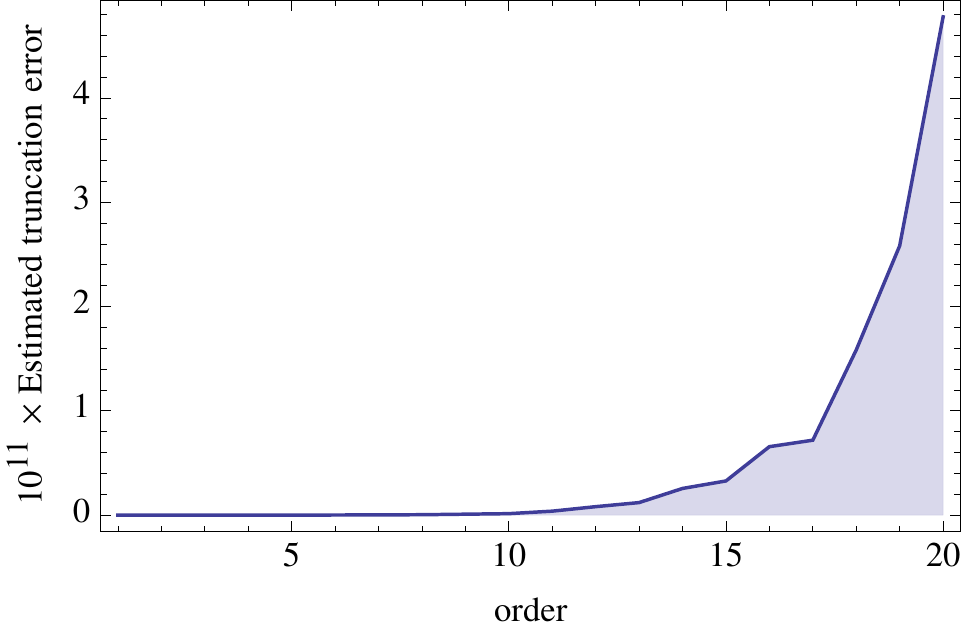}
\caption{Growth of the truncation errors due to the floating point arithmetic.}
\label{f:ratioRealComplex}
\end{figure}

The propagation of the truncation errors when using floating point arithmetic can be studied with the help of interval arithmetic \cite{CellettiChierchia1988,Jorba1999}. An alternative way of estimating these errors is as follows. On the one hand, the transformation from complex to real variables, which is exact when avoiding decimal expansions, will produce some residual complex terms due to the floating point arithmetic ---which, of course, must be neglected. The absolute value of the greatest of the coefficients affecting these residual terms is an indicator of the truncation errors accumulated in the computations. On the other hand, the size of the coefficients of the monomials generally grows with the consecutive higher orders of the perturbation solution. Then, the ratio between the greatest coefficient of the spurious, complex monomials and the greatest coefficient of the true, real monomials can be taken as an estimator of the truncation errors introduced by the computer's arithmetic. This is illustrated in Fig.~\ref{f:maxCoefW}, where we see that the growth rate of the complex residuals is higher than that of the coefficients of the real monomials, and Fig.~\ref{f:ratioRealComplex}, where we note the sharp growth of the truncation errors estimated with our criterion when passed the order 15th. 
\par

Besides, in view of the already mentioned rapid growth of the number of terms to be evaluated and because the dynamics about the libration points is generally highly unstable, thus making orbit propagation quite sensitive to the initial conditions, we propose values $\mathcal{O}(10^{-12})$ of this indicator as a practical limit for the validity of the analytical solution. From Fig.~\ref{f:ratioRealComplex}, this value would correspond to an order of the perturbation solution between, say, 10 and 16, yet these high orders are only needed for computing orbits far away from the libration points.
\par

\section{Performance of the analytical solution} \label{se:performance}

The reduction carried out by the normalization is more insightfully appreciated when using Lissajous canonical variables $(\ell,g,L,G)$ \cite{Deprit1991}. In theses variables
\begin{eqnarray} \label{I0}
I_0 &=& \mbox{$\frac{1}{2}$}L, \\ \label{I1}
I_1 &=& \omega{s}d\cos2g, \\ \label{I2}
I_2 &=& \omega{s}d\sin2g, \\ \label{I3}
I_3 &=& \mbox{$\frac{1}{2}$}G, 
\end{eqnarray}
with
\begin{equation} \label{sdLi}
s=\sqrt{\frac{L+G}{2\omega}}, \qquad d=\sqrt{\frac{L-G}{2\omega}},
\end{equation}
showing that the normalization removed the angle $\ell$, the conjugate variable to the momentum $L$. Then, the only variables of the reduced phase space are $g$ and $G$, and the orbits are ellipses whose size, shape, and orientation evolve slowly.
\par

The periodic orbits of the original space (the Hill problem) are computed analytically as follows. First, we compute the equilibrium in the Hopf variables representing the desired periodic orbit (planar or vertical Lyapunov orbits, Halo orbits, or periodic orbits of the bridge linking planar and vertical Lyapunov orbits). Then, $L$ and $G$ are trivially obtained from Eqs.~(\ref{I0}) and (\ref{I3}), respectively, while $g$ is computed unambiguously from Eqs.~(\ref{I1}) and (\ref{I2}). The choice of any particular value $\ell\in[0,2\pi)$ allows for the following computation of the complex variables
\begin{eqnarray} \label{v2Li}
 v &=& \sqrt{\frac{\omega}2}\left[(s-d)\cos{g}+\ii(d+s)\sin{g}\right](\cos\ell+\ii\sin\ell), \\ \label{w2Li}
 w &=& \sqrt{\frac{\omega}2}\left[(s-d)\sin{g}-\ii(d+s)\cos{g}\right](\cos\ell+\ii\sin\ell), \\ \label{vv2Li}
 V &=& \sqrt{\frac{\omega}2}\left[\ii(d-s)\cos{g}-(d+s)\sin{g}\right](\cos\ell-\ii\sin\ell), \\ \label{ww2Li}
 W &=& \sqrt{\frac{\omega}2}\left[\ii(d-s)\sin{g}+(d+s)\cos{g}\right](\cos\ell-\ii\sin\ell).
\end{eqnarray}
Next, the Lie transformation computed for achieving the reduction provides corresponding complex prime variables, from which the subindex 1 Cartesian variables are recovered using Eq.~(\ref{tocomplex}). Finally, Eq.~(\ref{x1tox}) ---with the matrix $A$ given by traditional choice in Eq.~(\ref{Agerard}) or any other choice from the family represented by Eq.~(\ref{matrixA}) with the constraints in Eq.~(\ref{constraint}) that could have been used alternatively--- will provide the initial conditions relative to the libration point. Repetition of the procedure for different values of $\ell$ will give the desired orbit without need of integrating these initial conditions.
\par

We explore the performance of the analytical solutions by comparing orbits predicted by different orders of the perturbation solution with their partner periodic orbits of the Hill problem computed numerically. We do the comparisons in three different scenarios. In the first one, we constrain the energy to values close to the energy of the libration point, a case in which only Lyapunov orbits exist. In the second case we explore higher energies, for which Halo orbits also exist but remain close to the libration point. Finally, we focus on the range of energy values in which, besides the Lyapunov and Halo orbits, the orbits of the two-lane bridge linking planar and vertical Lyapunov orbits exist. The later is a quite challenging case because of the large size of the orbits, and, until our knowledge, has never computed before analytically.

\subsection{$L=0.01$}

\begin{figure}[htb]
\centering \includegraphics[scale=0.85, angle=0]{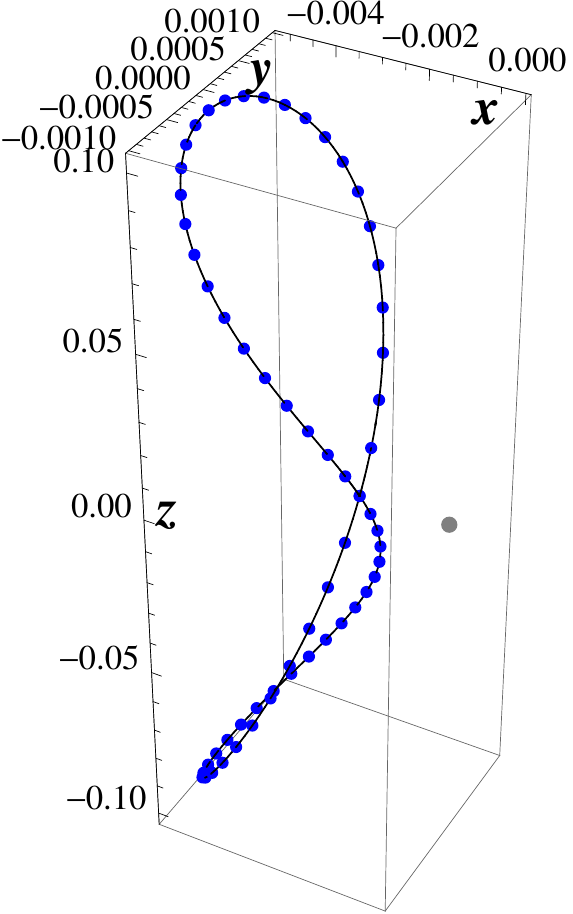} \qquad \includegraphics[scale=0.78, angle=0]{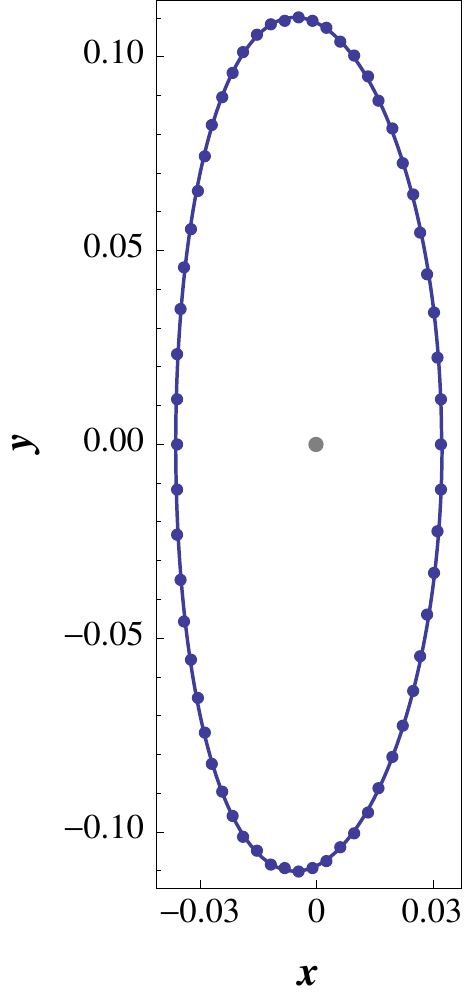}
\caption{Lyapunov orbits predicted by the analytical solution for $L=0.01$ (dots) superimposed to the corresponding, numerically integrated, periodic orbits (full line). In this and following plots, the gray dots represent the Lissajous point. Note that the $z$ axis of the left plot is in a different scale from that of the $x$ and $y$ axis.}
\label{f:Lyapunov01}
\end{figure}

For the vertical Lyapunov orbit, the order 4 of the perturbation solution is enough to mimic the true periodic orbit at the precision of the graphics, as shown in the left plot of Fig.~\ref{f:Lyapunov01}. However, when initial conditions provided by the analytical solution are propagated in the original, Hill problem dynamics, the periodicity error $\epsilon$, which is defined as 
\begin{equation}
\epsilon=\max|\beta_i(T)-\beta_i(0)|\qquad (i = 1,...,6),
\end{equation}
where $\beta_i$ stands for any of the coordinates in the original phase space, is only of the order of $10^{-6}$. The accuracy increases in a continuous way up to the order 11, for which the propagation in the original model of initial conditions taken from the analytical solution result in a periodicity error better than $10^{-10}$ after a period $T=3.146695654477$.
\par

In the case of the planar Lyapunov orbit, the order 3 of the perturbation solution suffices for suplying initial conditions that close the orbit at the precision of the graphics (right plot of Fig.~\ref{f:Lyapunov01}). However, $\epsilon$ is only of the order of $10^{-4}$. The periodicity error improves with higher orders of the solution and $\epsilon$ is of the order of $10^{-11}$ for the order 11. Only very slight improvements are achieved when the perturbation solution is truncated to higher orders, which become negligible further than the order 13th.

\subsection{$L=0.2$}

\begin{figure}[htbp]
\centering \includegraphics[scale=0.75, angle=0]{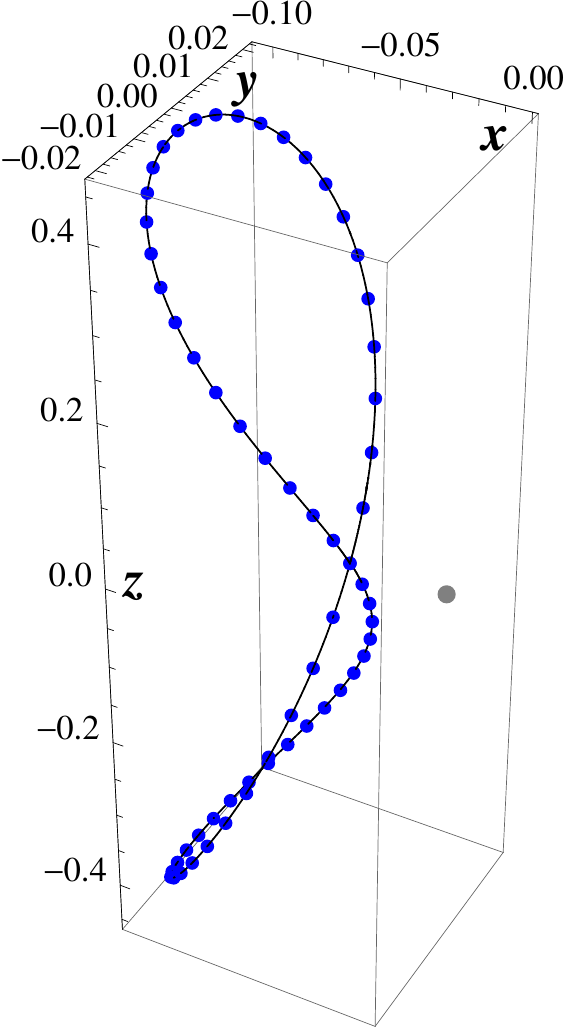} \qquad \includegraphics[scale=0.75, angle=0]{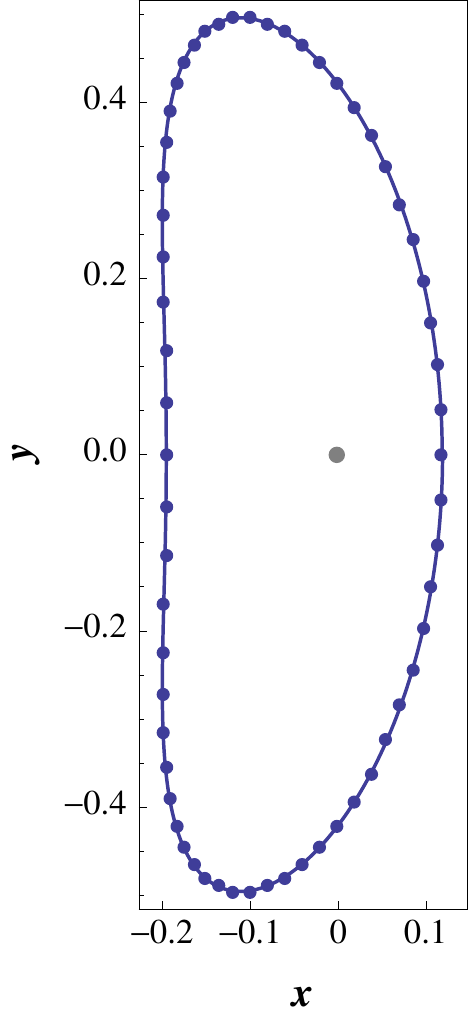}
\caption{Lyapunov vertical (left) and planar (right) orbits predicted by the analytical theory for $L=0.2$ (dots) superimposed to the corresponding, numerically integrated, periodic orbit (full line). Note the different scales of the axes of the left plot.}
\label{f:Lyapunov2}
\end{figure}

Due to the fact that the orbits are much larger in this case, the 7th order of the perturbation theory is required to match a vertical Lyapunov at the precision of the graphics (left plot of Fig.~\ref{f:Lyapunov2}), but the periodicity error is only $\epsilon=\mathcal{O}(10^{-3})$. This value improves for increasing orders of the perturbation theory up to the order 15th, where $\epsilon=\mathcal{O}(10^{-6})$, and does not improve with higher orders. For the planar Lyapunov orbit we needed to use the 9th order truncation of the theory to achieve initial conditions in the original problem leading to a periodicity error of the order of one thousandth, which is enough to close the orbit at the precision of the graphics, as shown in the right plot of Fig.~\ref{f:Lyapunov2}. Like in the case of the vertical Lyapunov orbit, the order 15th of the perturbation theory improves the periodicity up to $\epsilon=\mathcal{O}(10^{-6})$, but no further improvements are  found with higher orders of the solution.
\par

\begin{figure}[htbp]
\centering \includegraphics[scale=0.77, angle=0]{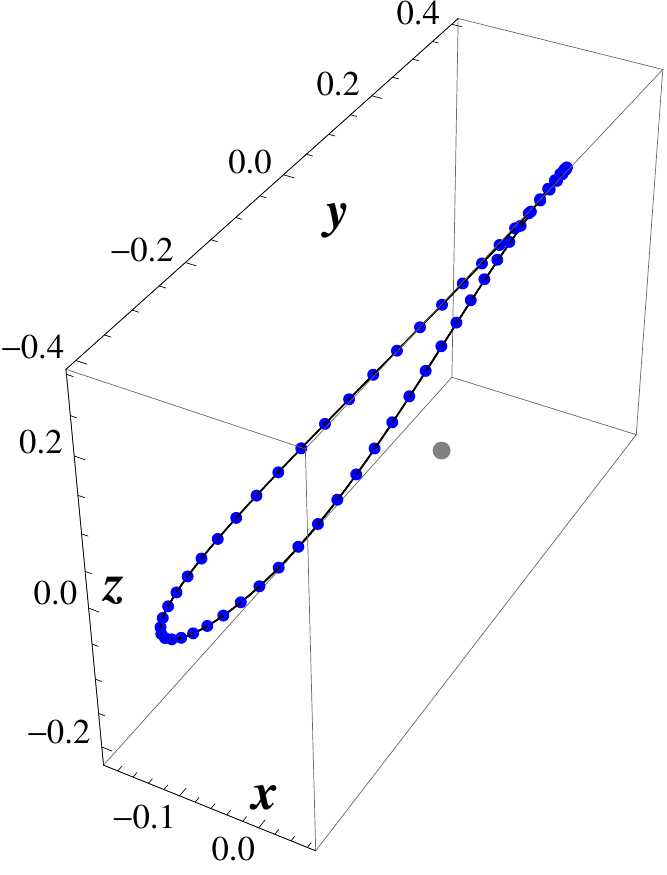} \, \includegraphics[scale=0.8, angle=0]{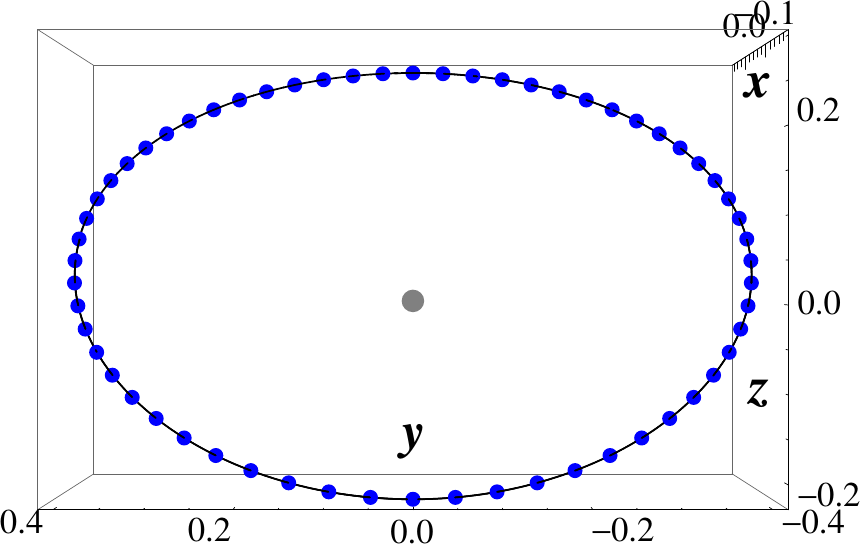}
\caption{Halo orbit predicted by the analytical solution for $L=0.2$ (dots) superimposed to the corresponding, numerically integrated, periodic orbit (full line).}
\label{f:Halo2}
\end{figure}

At this value of the energy, Halo orbits already bifurcated from the family of planar Lyapunov orbits. The analytical solution predicts them correctly starting from the ninth order truncation of the analytical solution, as shown in Fig.~\ref{f:Halo2}, where the periodicity error is $\epsilon=\mathcal{O}(10^{-3})$. Successive higher orders of the perturbation solution succeed in gradually improving periodicity, but only up to $\epsilon=\mathcal{O}(10^{-7})$, which happens with the order 17th of the analytical solution. Note that, while Lyapunov orbits always have the same initial conditions in the reduced phase space, as given by Eq.~(\ref{vyp}), the location of Halo orbits on the sphere is given by Eq.~(\ref{halo}) after solving Eq.~(\ref{haloRoot}), and, therefore, depends on the order of the perturbation solution used in each case.
\par

\subsection{$L=0.9$}
For this large value of $L$, the periodic orbits are so large that in different parts of the orbits the radius to the libration point falls out the convergence region of the Legendre polynomials expansion, which defines the validity of the model. In particular, the Halo orbit surrounds the central body, and its distance to the libration point is always longer than the Hill radius. Therefore, it is not expected that initial conditions obtained from the analytical theory can be improved by differential corrections to converge to a true Halo periodic orbit. 

Contrary to the Halo case, only parts of the other periodic orbits remain out the region where the Legendre polynomials expansion converges, and the analytical solutions succeeds in providing reasonable approximations to the true Lyapunov orbits. Thus, the order 13th of the perturbation solution is able to capture an approximation of the vertical Lyapunov orbit with periodicity $\epsilon=\mathcal{O}(10^{-4})$, yet no further improvements are found for higher orders of the perturbation solution. Analogously, initial conditions of a planar Lyapunov periodic orbit with periodicity of the order of $10^{-3}$ are obtained with the order 16th of the perturbation solution. In both cases the initial conditions are easily improved with differential corrections, leading to the true periodic orbits.
\par

\begin{figure}[htbp]
\centering \includegraphics[scale=0.65, angle=0]{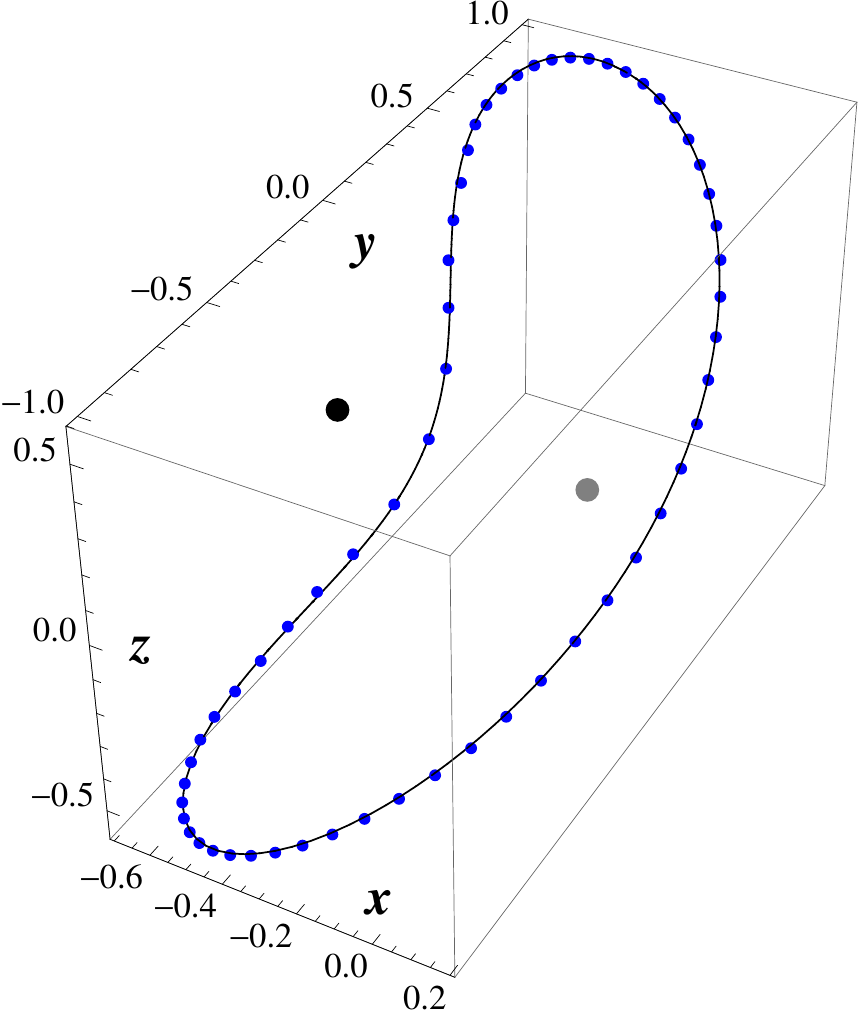}  \includegraphics[scale=0.7, angle=0]{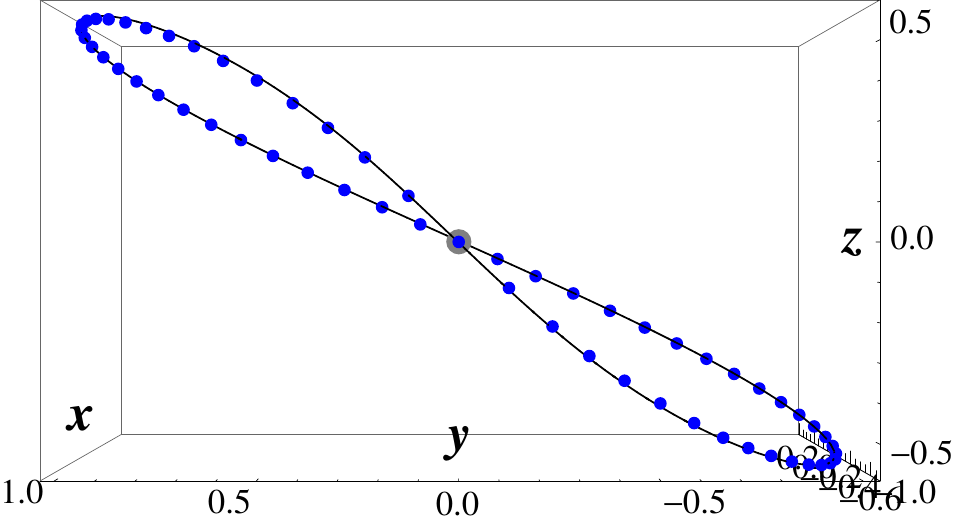}
\caption{Two different viewpoints of the orbit of the bridge family linking planar and vertical Lyapunov orbits predicted by the 14th order of the analytical solution when $L=0.9$ (dots) superimposed to the corresponding, numerically integrated, periodic orbit (full line). 
The gray dot is the libration point and the black dot is the primary.
}
\label{f:bridge9}
\end{figure}

Furthermore, the second bifurcation of the family of Lyapunov planar orbits has already happened at this value of $L$, and the perturbation solution is effective in capturing an orbit of the two-lane bridge that links planar and vertical Lyapunov orbits. Indeed, as shown in Fig.~\ref{f:bridge9}, a truncation to the order 14th of the perturbation solution provides an orbit with periodicity $\epsilon=10^{-2}$. While the periodicity error is not improved with higher orders of the analytical solution, initial conditions provide by the order 14th are easily improved by differential corrections to get the true periodic solution. 
\par

\section{Conclusions}

A higher order normalization of the Hill problem Hamiltonian centered at a libration point eases the computation of a single perturbation solution that captures the four main families of periodic orbits of the Hill problem originated from the libration points. Namely, the families of planar and vertical Lyapunov orbits, the family of Halo orbits, and the two-lane bridge of periodic orbits that connects both families of Lyapunov orbits. Planar and vertical Lyapunov orbits exist for all energies above the energy of the libration points, and, therefore, close to the libration points these kinds of orbits are accurately reproduced with the lower orders of the analytical solution. Still, higher orders of the solution are required if one wants to obtain vertical and planar Lyapunov orbits far away from the libration points within an acceptably accuracy. The initial orbits of the Halo family require, at least, the 5th order truncation of the perturbation solution, whereas the orbits of the two-lane bridge of periodic orbits that connect the families of planar and vertical Lyapunov orbits require, at least, a 14th order truncation of the analytical solution to obtain a reasonable approximation of the corresponding real periodic orbits, whose initial conditions can be improved by means of differential corrections to get the true periodic solution.
\par

The Hill problem has been chosen to study analytically the dynamics about the libration points because of its generality and simplicity. However, the procedures used in this research are general and can be analogously applied to the restricted three-body problem or variations of it including different perturbations.

\section*{Acknowledgemnts}

Partial support by the Spanish State Research Agency and the European Regional Development Fund under Project ESP2016-76585-R (AEI/ERDF, EU) is recognized. The first author (ML) was also supported by project ESP2013-41634-P of the same agencies.

\end{document}